\documentclass[12pt]{article}
\usepackage{mathrsfs}
\usepackage{amssymb}
\usepackage{t1enc}
\usepackage[latin1]{inputenc}
\usepackage[english]{babel}
\usepackage{amsmath,amsthm}
\usepackage{amsfonts}
\usepackage{latexsym}
\usepackage[dvips]{graphicx}
\usepackage{geometry,graphics,enumerate,tabularx,shapepar}
\usepackage[all,2cell,dvips]{xy} \UseAllTwocells \SilentMatrices
\date{}

\begin{document}
\title{Cluster-tilted algebras of type $D_n$$^\star$}
\author{Wenxu Ge, Hongbo Lv, Shunhua Zhang$^*$ \\
{\small Department of Mathematics,\ Shandong University,\ Jinan
250100, P. R. China} \\
{\small Dedicated to Professor Shaoxue Liu on the occasion of his
eightieth birthday}}

\maketitle

\begin{abstract}
Let $H$ be a hereditary algebra of Dynkin type $D_n$ over a field
$k$ and $\mathscr{C}_H$ be the cluster category of $H$. Assume that
$n\geq 5$ and that $T$ and $T'$ are tilting objects in
$\mathscr{C}_H$. We prove that the cluster-tilted algebra
$\Gamma=\mathrm{End}_{\mathscr{C}_H}(T)^{\rm op}$ is isomorphic to
$\Gamma'=\mathrm{End}_{\mathscr{C}_H}(T')^{\rm op}$ if and only if
$T=\tau^iT'$ or $T=\sigma\tau^jT'$ for some integers $i$ and $j$,
where $\tau$ is the Auslander-Reiten translation and $\sigma$ is the
automorphism of $\mathscr{C}_H$ defined in section 4.
\end{abstract}


\maketitle


\footnote{ $^\star$ Supported by the NSF of China (Grant No.
10771112) and of the NSF of Shandong province (Grant No. Y2008A05).}

\footnote{ $^{*}$ Corresponding author.}

\footnote{ {\it Email addresses}: \  gewenxu1982@163.com(W.Ge), \
 lvhongbo356@163.com(H.Lv), \\  shzhang@sdu.edu.cn(S.Zhang)}

\section {Introduction}
Cluster algebras were introduced by S. Fomin and A. Zelevinsky [FZ1,
FZ2]. Cluster category was defined in [BMRRT] as a means for a
better understanding of cluster algebras of Fomin and Zelevinsky,
which is an orbit category $\mathscr{C}_H=\mathscr{D}^b({\rm mod}\
H)/\tau^{-1}[1]$ of the derived category of a hereditary algebra $H$
and is a triangulated category, see also [CCS1] for the cluster
categories of type $A_n$. Now cluster category becomes a successful
model for acyclic cluster algebras.

\vskip 0.2in

The cluster-tilted algebra is an algebra of the form
$\Gamma=\mathrm{End}_{\mathscr{C}_H}(T)^{op}$, where $T$ is a
tilting object in $\mathscr{C}_H$. The module category of finitely
generated $\Gamma$-modules mod-$\Gamma$ was explicitly described in
[BMR1], and the cluster-tilted algebra is considered as a path
algebra of a quiver with relations in [BMR2], and it was proved that
a (basic) cluster-tilted algebra of finite type is uniquely
determined by its quiver. Moreover, the relations were explicitly
described in [CCS2] for a cluster-tilted algebra of finite type.

\vskip 0.2in

A geometric realization of cluster categories of type $A_n$ was
given in [CCS1], and a calculating formula for the number of
non-isomorphic cluster-tilted algebras of type $A_n$ is given in
[Tor]. R.Schiffler gave a geometric realization for cluster
categories of type $D_n$ in [Sch]. In this paper, we use the
geometric realization of [Sch] to prove that two cluster-tilted
algebras $\mathrm{End}_{\mathscr{C}_H}(T)^{\rm op}$ and
$\mathrm{End}_{\mathscr{C}_H}(T')^{\rm op}$ are isomorphic if and
only if $T=\tau^iT'$ or $T=\sigma\tau^jT'$ for some integers $i$ and
$j$, where $\tau$ is the Auslander-Reiten translation and $\sigma$
is an automorphism of $\mathscr{C}_H$, see section 4.

\vskip 0.2in

This paper is arranged as follows. In section 2, we  collect
necessary definitions and basic facts needed for our research. In
section 3, we give an explicit description of the equivalence class
of triangulations of the category of tagged edges of the punctured
polygon ${\bf P}_n$. As an application, we obtain the mutation
classes of quivers of type $D_n$, and deduce all the quivers of
cluster-tilted algebras  of type $D_n$, and moreover, we also give
an explicit description for the relations which is consistent with
[CCS2]. In section 4, we prove that if  $T$ and $T'$ are basic
tilting objects in the cluster category  $\mathscr{C}_H$ of type
$D_n$, then the cluster-tilted algebras
$\mathrm{End}_{\mathscr{C}_H}(T)^{\rm op}$ and
$\mathrm{End}_{\mathscr{C}_H}(T')^{\rm op}$ are isomorphic if and
only if $T =\tau^i T'$ or $T=\sigma \tau^j T'$ for some integers $i$
and $j$.

\vskip 0.2in

Throughout this paper, we fix a field $k$, and denote by  $H$ a
hereditary $k$-algebra of type $D_n$. We follow the standard
terminology and notation used in the representation theory of
algebras, see [ARS], [ASS], [Hap] and [Rin].

\section {Preliminaries}

\vskip 0.2in

 Let $n$ be an integer with $n\geq 4$ and $\mathbf{P}_n$ be
 the punctured polygon with one puncture in its center. We give an
 example in Figure 1 for $n=8$.

 \begin{figure}[htbp]
\begin{center}
\includegraphics[width=3cm,height=3cm]{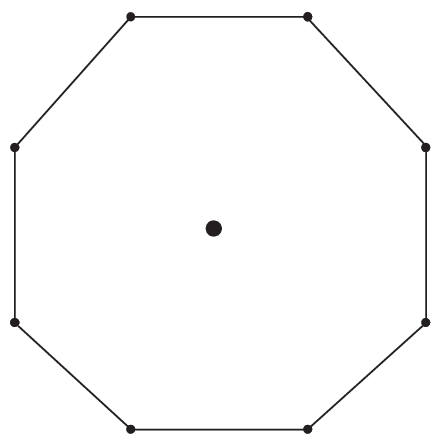} \\
{Figure 1}
\end{center}
\end{figure}

 We recall some definitions and results from [Sch] which will be needed
 in our further research.

\vskip 0.2in

 Let  $a\neq b$ be vertices on the boundary of $\mathbf{P}_n$ and
$\delta_{a,b}$ be the path along the boundary from $a$ to $b$ in
counterclockwise direction which does not run through the same point
twice and $|\delta_{a,b}|$ be the number of vertices on the path
$\delta_{a,b}$ (including $a$ and $b$). Two vertices $a$ and $b$ are
called neighbors if $|\delta_{a,b}|=2$ or $|\delta_{a,b}|=n$, and
$b$ is the counterclockwise neighbor of $a$ if $|\delta_{a,b}|=2$.
Note that $\delta_{a,a}$ is the path along the boundary from $a$ to
$a$ in counterclockwise direction which goes around the polygon
exactly once and such that $a$ is the only point through which
$\delta_{a,a}$ runs twice, and then $|\delta_{a,a}|=n+1$.

\vskip 0.2in

An edge of $\mathbf{P}_n$ is a triple $(a,\alpha,b)$ where $a$ and
$b$ are vertices on the boundary and $\alpha$ is a path from $a$ to
$b$ such that  $\alpha$ does not cross itself and $\alpha$ lies in
the interior of $\mathbf{P}_n$ (except for its starting point $a$
and its endpoint $b$) and $\alpha$ is homotopic to $\delta_{a,b}$
with $|\delta_{a,b}| \geq3$.

\vskip 0.2in

Two edges $(a,\alpha,b)$, $(c,\beta,d)$ are equivalent if $a=c$,
$b=d$, and $\alpha$ is homotopic to $\beta$. We denote by $E$ the
set of equivalence class of edges. Since an element of $E$ is
uniquely determined by an ordered pair of vertices $(a,b)$, we will
therefore use the notation $M_{a,b}$ for the equivalence class of
edges $(a,\alpha,b)$ in $E$.

\vskip 0.2in

Let $E'=\{M_{a,b}^\epsilon \ |\  M_{a,b} \in E,\  \epsilon = \pm 1\
{\rm and}\ \ \epsilon = 1\  {\rm if} \ a\neq b\}$ be the set of
tagged edges. If $a \neq b$, we will write $M_{a,b}$ instead of
$M_{a,b}^1$, and say that $M_{a,b}$ is a tagged edge with the length
of $|\delta_{a,b}|$. Note that there are exactly two tagged edges
$M_{a,a}^{-1}$ and $M_{a,a}^{1}$ for every vertex $a$ on the
boundary of $\mathbf{P}_n$, and the tagged edges $M_{a,a}^\epsilon$
will not be represented as loops around the puncture but as lines
from the vertex $a$ to the puncture. If $\epsilon = -1$, the line
will have a tag on it and if $\epsilon = 1$, there will be no tag.

\vskip 0.2in

We use the convention that $M_{a,b}$ is a tagged edge with  length
of $|\delta_{a,b}|$ if $a \neq b$, and that $M_{a,a}^\epsilon$ is a
tagged edge with length of $1$. We say that $M_{a,b}$ is \emph{close
to the border} if its length is three, and that $M_{a,b}$ is
\emph{connected} if its length is at least four. For every vertex
$a$ on the boundary of $\mathbf{P}_n$,  we say that
$M_{a,a}^\epsilon$ is \emph{degenerate}.

\vskip 0.2in

Let $M=M_{a,b}^\epsilon$ and $N=M_{c,d}^{\epsilon'}$ be in $E'$. The
crossing number $e(M,N)$ of $M$ and $N$ is the minimal number of
intersection of $M_{a,b}^\epsilon$ and $M_{c,d}^{\epsilon'}$ in
$\Delta^o$, where $\Delta^o$ is the interior of the punctured
polygon $\mathbf{P}_n$.  A {\it triangulation } of the punctured
polygon $\mathbf{P}_n$ is a maximal set of non-crossing tagged
edges.

\vskip 0.2in

{\bf Lemma 2.1.}$^{[{\rm Sch}]}$ \ (1)\ {\it Any triangulation of
  $\mathbf{P}_n$ has $n$ elements.}

 \vskip 0.1in

 (2)\  {\it $e(M,M)=0$ for $M\in E'$ and
\ $e(M_{a,a}^{\epsilon}, M_{b,b}^{\epsilon'}) =\left\{
          \begin{array}{ll}
          {1 \ \ {\rm if} \ a\neq b\
{\rm and} \  \epsilon\neq \epsilon'}\\
           {0 \ \ {\rm otherwise},}
          \end{array}
        \right.$}

\vskip 0.2in

According to [Sch], the translation $\tau$ is a bijection $\tau$:
$E'\rightarrow E'$ defined as the following.  Let
$M_{a,b}^\epsilon\in E'$, and let $a'$ (respectively $b'$) be the
clockwise neighbor of $a$ (respectively $b$). Then  $\tau M_{a,b}=
M_{a',b'}$ if $a\neq b$, and $\tau M_{a,a}^\epsilon=
M_{a',a'}^{-\epsilon}$ if $a=b$ with $\epsilon=\pm 1$ .

\vskip 0.2in

Let $\mathcal{C}$ be the $k$-linear additive category whose objects
are direct sums of tagged edges in $E'$, and the morphism from a
tagged edge $M$ to a tagged edge $N$ is a quotient of the vector
space over $k$ spanned by sequences of elementary moves starting at
$M$ and ending at $N$ see 3.6 in [Sch] for details.

\vskip 0.2in

Let $\mathscr{C}_H$ be the cluster category of $H$, where $H$ is a
hereditary $k$-algebra of type $D_n$.  It is well known that
$\mathscr{C}_H=\mathscr{D}^b({\rm mod}\ H)/\tau^{-1}[1]$ is a
triangulated category. A basic {\it tilting object} in the cluster
category $\mathscr{C}_H$ is an object $T$ with $n$ non-isomorphic
indecomposable direct summands such that ${\rm
Ext}^1_{\mathscr{C}_H}(T,T)=0$.

\vskip 0.2in

The following lemma is proved in [Sch] and it will be used later.

\vskip 0.2in

{\bf Lemma 2.2.} \  {\it There is an equivalence of categories
$\varphi$ between the category of tagged edges $\mathcal{C}$ and the
cluster category $\mathscr{C}_H$ of type $D_n$, moreover,
$\mathcal{C}$ is a triangulated category with the shift functor
$[1]=\tau$.}

\vskip 0.2in

According to [Sch], We can define the $\mathrm{Ext}^1$ of two
objects $M,N\in {\rm ind }\ \mathcal{C}$ as
$\mathrm{Ext}_\mathcal{C}^1(M,N)=\mathrm{Hom}_\mathcal{C}(M,\tau
N)$.

\vskip 0.2in

{\bf Lemma 2.3.} \  {\it  Let $M, N\in {\rm ind}\ \mathcal{C} $.
Then ${\rm dim}\ \mathrm{Ext}_\mathcal{C}^1(M,N)= e(M,N)$.}

\vskip 0.2in

 Let $\mathcal {T}_n$ be the set of all triangulations of the punctured polygon
$\mathbf{P}_n$. According to Lemma 2.2 and Lemma 2.3, we know that
$\varphi(\Delta)$ is a basic tilting object in $\mathscr{C}_H$ for
any triangulation $\Delta\in \mathcal {T}_n $. We denote by
$Q_\Delta$ the quiver of the cluster-tilted algebra
$\mathrm{End}_{\mathscr{C}_H}(\varphi(\Delta))^{\rm op}$.  It is
well known that the vertices of $Q_\Delta$ is the isomorphism
classes of simple
$\mathrm{End}_{\mathscr{C}_H}(\varphi(\Delta))^{\rm op}$ modules,
and the number of the arrows from simple modules $S_1$ to $S_2$ is
${\rm dim}_k {\rm Ext}^1(S_1,S_2)$. It is easy to see that all
quivers  of cluster-tilted algebras of type $D_n$ can be obtained in
this way.

\vskip 0.2in

We define a mutation of a triangulation at a given tagged edge by
replacing this tagged edge with another tagged edge described as 6.2
in [Sch]. Let $Q_\Delta$ be the quiver corresponding to a
triangulation $\Delta\in \mathcal {T}_n$. Then the mutation of
$Q_\Delta$ at the vertex $i$ corresponds to the mutation of $\Delta$
at the tagged edge corresponding to $i$. It is easy to see that each
triangulation induces a quiver of a cluster-tilted algebra, and that
every quiver of cluster-tilted algebras of type $D_n$ can be
assigned to at least one triangulation.

\vskip 0.2in

Let $\mathcal {M}_n^D$ be the mutation class of type $D_n$. We have
a function $\gamma$: $\mathcal {T}_n\rightarrow \mathcal {M}_n^D$
defined by $\gamma(\Delta)=Q_\Delta$ for any triangulation $\Delta$
in $\mathcal {T}_n$, it is easy to see that $\gamma$ is surjective.
Let $\mathcal{H}_n$ be the isomorphism class of cluster-tilted
algebras of type $D_n$. Note that any cluster-tilted algebra of type
$D_n$ is uniquely determined by its quiver up to isomorphism and
then we have a function $\gamma'$: $\mathcal{T}_n\rightarrow
\mathcal{H}_n$ which is induced by $\gamma$. Note that $\gamma'$ is
also surjective.

\vskip 0.2in

{\bf Definition 2.4.} {\it  Let $\mathcal{C}$ be the category of
tagged edges of $\mathbf{P}_n$. We define a function $\sigma$:
$\mathcal{C}\rightarrow\mathcal{C}$  by $\sigma(M_{a,b})=M_{a,b}$ if
$a\neq b$, and $\sigma(M_{a,b}^\epsilon)=M_{a,b}^{-\epsilon}$ if
$a=b$ with $\epsilon=\pm 1$.}

\vskip 0.2in

Note that $\sigma$ is an automorphism and $\sigma^2=\mathbf{1}$. For
any triangulation $\Delta\in \mathcal {T}_n$ we have that
$\gamma(\Delta)=\gamma(\sigma\Delta)$ and $\tau\sigma=\sigma\tau$.

\vskip 0.2in

We define an equivalence relation $\sim$ on $\mathcal {T}_n$ as the
follows: $\Delta\sim\Delta'$ if $\Delta=\tau^i \Delta'$ or
$\Delta=\sigma \tau^j \Delta'$ for some integers $i$ and $j$.  We
denote by  $\widetilde{\mathcal {T}_n}$ the equivalence class of
$\mathcal {T}_n$, and we have a  function $\tilde{\gamma}$:
$\widetilde{\mathcal {T}_n}\rightarrow \mathcal{M}^D_n$ induced by
$\gamma$.  Note that $\tilde{\gamma}$ is surjective and well
defined, and that $Q_\Delta=Q_{\Delta'}$ in $\mathcal{M}^D_n$
whenever $\Delta\sim\Delta'$.  In section 4, we will prove that
$\tilde{\gamma}$ is also injective if $n\geq 5$.

\vskip 0.2in

\section {Quivers of cluster-tilted algebras of type $D_n$}

In this section, we will give a classification for the elements in
$\widetilde{\mathcal {T}_n}$, and then  deduce an explicit
description for quivers of cluster-tilted algebras of type $D_n$.
That is, we give a new proof for the mutation class of quivers of
type $D_n$ obtained in [Vat]. We should mention that our proof is
started from the classification for the elements in
$\widetilde{\mathcal {T}_n}$ which seems to have independent
interest.

\vskip 0.2in

Let $\Delta$ be a triangulation in $\mathcal {T}_n$ and $Q_\Delta$
be the quiver of cluster-tilted algebra
$\mathrm{End}_{\mathscr{C}_H}(\varphi(\Delta))^{\rm op}$. Let $M$ be
a tagged edge in  $\Delta$. We denote by $V_M$ the vertex of
$Q_\Delta$ corresponding to $M$.  We say that $V_M$ is
\emph{connected} if the corresponding tagged edge $M$ in $\Delta$ is
connected. Two tagged edges $M$ and $N$ are adjacent in $\Delta$ if
the corresponding vertices $V_M$ and $V_N$ are adjacent by an arrow
in the corresponding quiver $Q_\Delta$.

\vskip 0.2in

Let $\Delta$ be a triangulation and $M=M_{a,b}\in\Delta$ be
connected. Then $M_{a,b}$ divides the punctured polygon
$\mathbf{P}_n$ into two parts,  one part $\mathbf{P'}$ without the
puncture in its interior and the other part $\mathbf{P''}$ with the
puncture. $M_{a,b}$ also divides the triangulation $\Delta$ into two
parts,  one part $\Delta_1$ is in $\mathbf{P'}$ and the other part
$\Delta_2$ is in $\mathbf{P''}$, and
$\Delta=\Delta_1\bigcup\Delta_2\bigcup\{M_{a,b}\}$. We have the
following lemma.

\vskip 0.2in

{\bf Lemma 3.1.}\ {\it  Let $\Delta_1$ and $\Delta_2$ be the
partition defined as above. Assume that $M_1\in\Delta_1$ and
$M_2\in\Delta_2$, and let $V_{M_1}$ and $V_{M_2}$ be the
corresponding vertices in the corresponding quiver $Q_\Delta$.  Then
$V_{M_1}$ is not adjacent to
 $V_{M_2}$ in  $Q_\Delta$, and  any path between $V_{M_1}$ and
$V_{M_2}$ must run through the vertex $V_M$.}

\vskip 0.1in

{\bf Proof.}  We assume by contrary that there is an arrow from
$V_{M_1}$ to $V_{M_2}$. By using Lemma 2.3, we have that
$e(M_1,\tau^{-1}M_2)=\mathrm{Ext}_\mathcal{C}^1(M_1,\tau^{-1}M_2)
=\mathrm{Hom}_\mathcal{C}(M_1,M_2)\neq 0$. That is,  $M_1$ and
$\tau^{-1}M_2$ intersect in the interior of the puncture polygon
$\mathbf{P}_n$. It is easy to see that  $M$ and $\tau^{-1}M_2$
intersect in the interior of the punctured polygon, that is
$e(M,\tau^{-1}M_2)=\mathrm{Hom}_\mathcal{C}(M,M_2)\neq 0$. By the
same argument we have that $\mathrm{Hom}_\mathcal{C}(M_1,M)\neq 0$.
Therefore, we have a path $V_{M_1}\rightarrow V_1 \rightarrow V_2
\rightarrow, ..., \rightarrow V_i \rightarrow V_M$ from $V_{M_1}$ to
$V_M$ and a path $V_M \rightarrow V'_1 \rightarrow, ..., V'_j
\rightarrow V_{M_2}$ from $V_M$ to $V_{M_2}$ in $Q_\Delta$.  Hence
we have a subquiver of $Q_\Delta$ as Figure 2.

\begin{figure}[htbp]
\begin{center}
\includegraphics[width=5cm,height=5cm]{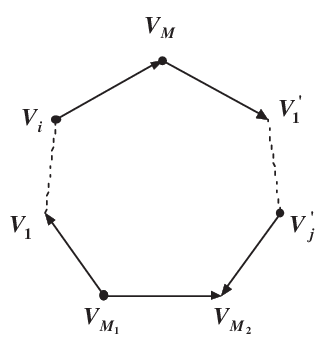} \\
{Figure 2}
\end{center}
\end{figure}

Perform the following mutations  $\mu_{V_1}$, $\mu_{V_2}$,... ,
$\mu_{V_i}$, $\mu_{V_M}$, $\mu_{V'_1}$,... , $\mu_{V'_j}$ one by
one, we then obtain a multiple arrows from $V_{M_1}$ to $V_{M_2}$,
which is a contradiction. This completes the proof.  $\Box$

\vskip 0.2in

{\bf  Remark.}\ We denote by $m$ the number of vertices of
$\mathbf{P'}$ and write $\mathbf{P'}_m$ instead of $\mathbf{P'}$. We
obtain a subquiver $Q'_m$ of $Q_\Delta$ by factoring out all
vertices corresponding to $\Delta''$. According to [CCS1], $Q'_m$ is
a quiver of some cluster-tilted algebra of type $A_m$.

\vskip 0.2in

{\bf Lemma 3.2.}\ {\it Let $M$ be a connected tagged edge in a
triangulation $\Delta$. Let $Q'_m$ be a subquiver defined as above.
Then the number of vertices of $Q'_m$ which are adjacent to $V_M$ is
either one or two. If this number is two, then $V_M$ lies on an
oriented cycle of length three.}

\vskip 0.1in

{\bf Proof.} \  Let $V_N$ be a vertex in the quiver $Q'_m$ which is
adjacent to $V_M$.  According to Lemma 3.1, we only need to consider
the following three cases in Figure 3.

\begin{figure}[htbp]
\begin{center}
\includegraphics[width=10cm,height=3cm]{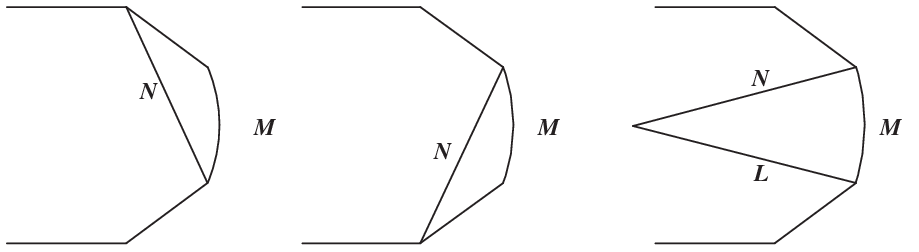} \\
{Figure 3}
\end{center}
\end{figure}

In the first two cases, the vertex which is adjacent to $V_M$ is
only $V_N$, and in the third case,  $V_M$ lies on the cycle $V_M
\rightarrow V_{L} \rightarrow V_N \rightarrow V_M$. $\Box$

\vskip 0.2in

Let $\mathcal{M}_k^A$ be the mutation class of type $A_k$. The union
of all $\mathcal{M}_k^A$ with $k\geq 1$ will be denoted by
$\mathcal{M}^A$. If $Q$ is a quiver in $\mathcal{M}^A$,  according
to [BV], we have the following facts.

\vskip 0.1in

         (1) \ The length of all oriented cycles is three.

\vskip 0.1in

         (2) \ The number of vertices which is adjacent to a fixed vertex is at most
         four.

\vskip 0.1in

         (3)\  If the number is four, then two of its adjacent arrows belong to one
         3-cycle, and the other two belong to another 3-cycle.

\vskip 0.1in

         (4)\  If the number is three, then two of its adjacent arrows belong to one
         3-cycle, and the third arrow does not belong to any
         3-cycle.

\vskip 0.2in

Let $H_A$ be a cluster-tilted algebra of type $A_n$. Then $H_A=kQ/I$
with $Q$ belonging to $\mathcal{M}_n^A$, and the relation ideal $I$
is generated by the set of all length 2 paths of all 3-cycles in
$Q$. We denote the relation set by $\Sigma_Q$, see [BMR2, BV] for
details.

\vskip 0.2in

Now we want to classify the elements in $\widetilde{\mathcal
{T}_n}=\mathcal{T}_n / \sim$,  where $\mathcal{T}_n $ is the set of
all triangulations of the punctured polygon $\mathbf{P}_n$, and
$\Delta\sim\Delta'$ if and only if $\Delta=\tau^i \Delta'$ or
$\Delta=\sigma \tau^j \Delta'$ for some integers $i$ and $j$. By
abuse of language, we also denote by $\Delta$ the equivalence class
determined by itself. Note that we have a surjective function
$\tilde{\gamma}$: $\widetilde{\mathcal {T}_n}\rightarrow
\mathcal{M}_n^D$, and we will prove that  $\tilde{\gamma}$ is also
injective in next section if $n\geq 5$.  In particular, we obtain a
classification for $\mathcal{M}_n^D$.

\vskip 0.2in

Two degenerate tagged edges $M_{a,a}^1$ and $M_{a,a}^{-1}$ are
called a {\it double edges}, and two degenerate tagged edges
$M_{a,a}^{\epsilon}$ and $M_{b,b}^{\epsilon}$ are called a {\it
pairing edges} if $a$ and $b$ are neighbors on the punctured polygon
$\mathbf{P}_n$. Let $\Delta$ be an element in $\widetilde{\mathcal
{T}_n}$ and $\tilde{\gamma}(\Delta)=Q_{\Delta}$ be the corresponding
element in $\mathcal{M}^D_n$. According to the definition of
triangulation, we know that $\Delta$ has at least two degenerate
tagged edges.

\vskip 0.2in

{\bf Lemma 3.3.}\  {\it Let $\Delta$ be an element in
$\widetilde{\mathcal {T}_n}$, and let $M_{a,a}^1$ and $M_{b,b}^1$ be
two degenerate tagged edges in $\Delta$ and $b$ be a
counterclockwise vertex of $a$. If there is no other degenerate
tagged edge between $M_{a,a}^1$ and $M_{b,b}^1$ in $\Delta$, then
$M_{a,b}$ is a tagged edge in $\Delta$.}

\vskip 0.1in

{\bf Proof.} \   Assume by contrary that the tagged edge $M_{a,b}$
is not in $\Delta$. Then there exists a non-degenerate tagged edge
$M_{c,d}$ in $\Delta$ such that $e(M_{c,d},M_{a,b})\neq 0$. It is
easy to see that $e(M_{c,d},M_{a,a}^1)\neq 0$ or
$e(M_{c,d},M_{b,b}^1)\neq 0$ which contradicts with the definition
of $\Delta$. This completes the proof. \hfill$\Box$

\vskip 0.2in

{\bf Proposition 3.4.}\ {\it  Let $\Delta$ be an element in
$\widetilde{\mathcal {T}_n}$.}

\vskip 0.1in

 (1)\ {\it If $M_{a,b}$ is a tagged edge in  $\Delta$ with
$|\delta_{a,b}|=n$, then $\Delta$  has exactly two  degenerate
tagged edges forming a double or a pairing.}

\vskip 0.1in

(2)\  {\it Assume that $\Delta$ has exactly two degenerate tagged
edges forming a double. We denote the double by $M_{a,a}^1$ and
$M_{a,a}^{-1}$, If $\Delta$ has no tagged edge of length $n$, then
there exists a counterclockwise vertex $b$ ($b\neq a$) on the
boundary of $\mathbf{P}_n$ such that $M_{a,b}$ and $M_{b,a}$ belong
to $\Delta$.}

\vskip 0.1in

(3)\ {\it Assume that $\Delta$ has exactly two degenerate tagged
edges which are not a double. If $\Delta$ has no tagged edge of
length $n$, then the two degenerate tagged edges is not a pairing.}

 \vskip 0.1in

(4)\ {\it Assume that $\Delta$ has more than three degenerate tagged
edges. Then $\Delta$ has no tagged edge of length $n$.}

\vskip 0.1in

{\bf Proof.} \  (1) It follows from Lemma 2.1.

\vskip 0.1in

(2)\ Assume that $b$ is  the counterclockwise neighbor of $a$ on the
punctured polygon $\mathbf{P}_n$  such that $M_{a,b}$ is a tagged
edge in $\Delta$ with maximal length,  We want to prove that
$M_{b,a}$ is a tagged edge in $\Delta$.  Note that $b$ is not a
neighbor of $a$ since $\Delta$ has no tagged edge of length $n$. If
$M_{b,a}$ is not in $\Delta$.  Then there exists a non-degenerate
tagged edge $M_{c,d}$ in $\Delta$ such that $e(M_{c,d},M_{b,a})\neq
0$. It is easy to see that $e(M_{c,d},M_{a,b}^1)\neq 0$ which
contradicts with that $M_{a,b}\in \Delta$. Therefore, $M_{a,b}$ and
$M_{b,a}$ are tagged edges of $\Delta$.

\vskip 0.1in

(3) It follows from Lemma 3.3.

\vskip 0.1in

(4)\ It follows from the definition of triangulation.   \hfill$\Box$

\vskip 0.2in

Let $\Delta$ be an element in $\widetilde{\mathcal {T}_n}$, and
$Q_{\Delta}= \tilde{\gamma}(\Delta)$ be the corresponding quiver in
$\mathcal {M}_n^D$. We denote by $H_\Delta$ the cluster-tilted
algebra corresponding to  $Q_{\Delta}$. According to Proposition
3.4, we can divide elements in $\widetilde{\mathcal {T}_n}$ into the
following four types, according to Lemma 2.2, we also obtain a
similar classification for $\mathcal {M}_n^D$ as in [Vat].

\vskip 0.2in

\textbf{Type 1}: \ {\it  $\Delta$ has a tagged edge $M=M_{a,b}$ with
$|\delta_{a,b}|=n$ as in Figure 4.}

\begin{figure}[htbp]
\begin{center}
\includegraphics[width=8cm,height=5cm]{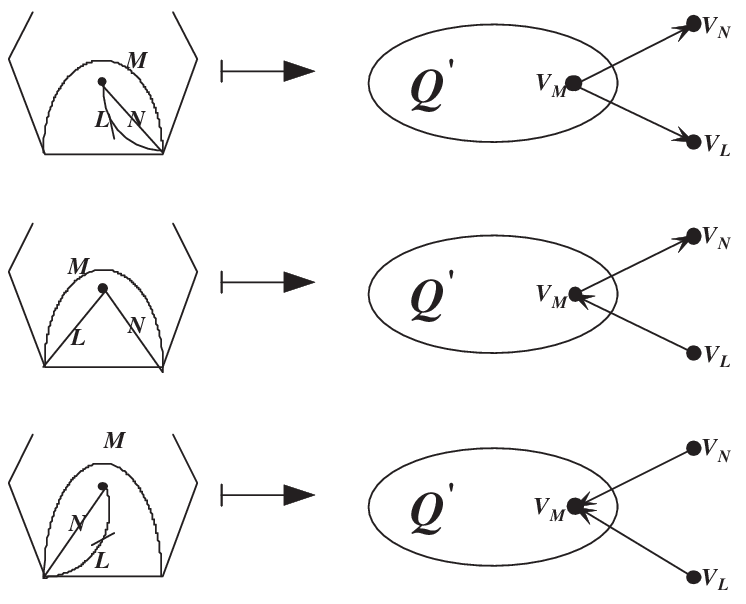} \\
{Figure 4}
\end{center}
\end{figure}

In this case, we write $\Delta$ and $Q_{\Delta}$ as in Figure 4.
Note that $Q'=Q_\Delta \setminus \{ V_N,\,V_{L}\}$ is in
$\mathcal{M}^A$ and $V_M$ is connected. The cluster-tilted algebra
$H_\Delta=kQ_\Delta/I$, the relation ideal $I$ is generated by the
set $\Sigma_{Q'}$, where $\Sigma_{Q'}$ is all the paths of length 2
which are part of a 3-cycle in $Q'$.

\vskip 0.2in

\textbf{Type 2}: {\it $\Delta$ has no tagged edge of length $n$ and
it has exactly two degenerate edges forming a double.}

\vskip 0.1in

In this case,  we may denote the double by $M_{a,a}^1$ and
$M_{a,a}^{-1}$. According to Proposition 3.4 (2), we know that there
exists a vertex $b$ ($b\neq a$) such that $M_{a,b}$ and $M_{b,a}$
belong to $\Delta$. Moreover, $M_{a,b}$ and $M_{b,a}$ divide the
triangulation $\Delta$ to three parts, $\Delta_1$, $\Delta_2$ and
$\{M_{a,b}, M_{b,a}, M_{a,a}^1, M_{a,a}^{-1}\}$.

\begin{figure}[htbp]
\begin{center}
\includegraphics[width=12cm,height=3cm]{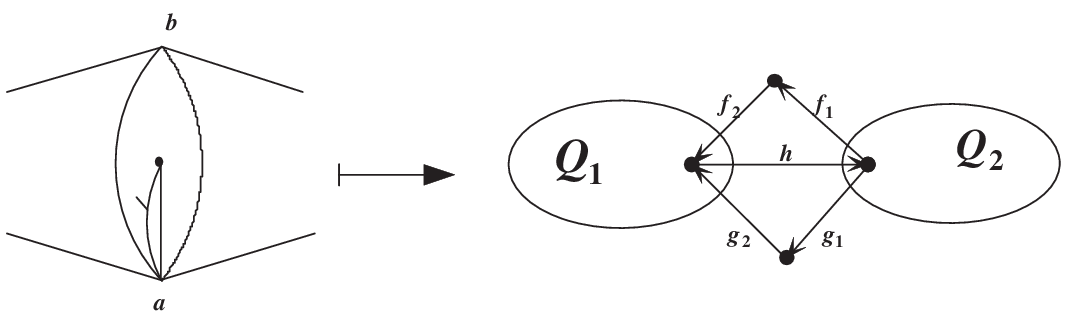} \\
{Figure 5}
\end{center}
\end{figure}

As in Figure 5 , $Q_1$ and $Q_2$ are in $\mathcal{M}^A$. The
cluster-tilted algebra $H_\Delta$ is the quotient algebra
$kQ_\Delta/I$, where the relation ideal $I$ is generated by the set
$\Sigma_{Q_1}\bigcup \Sigma_{Q_2} \bigcup \{f_1 f_2 - g_1 g_2, h
f_1, f_2 h, h g_1, g_2 h\}$

\vskip 0.2in

\textbf{Type 3}: {\it $\Delta$ has exactly two degenerate tagged
edges which are not a double (see Figure 6) and $\Delta$ has no
tagged edge of length $n$.}

\begin{figure}[htbp]
\begin{center}
\includegraphics[width=12cm,height=3cm]{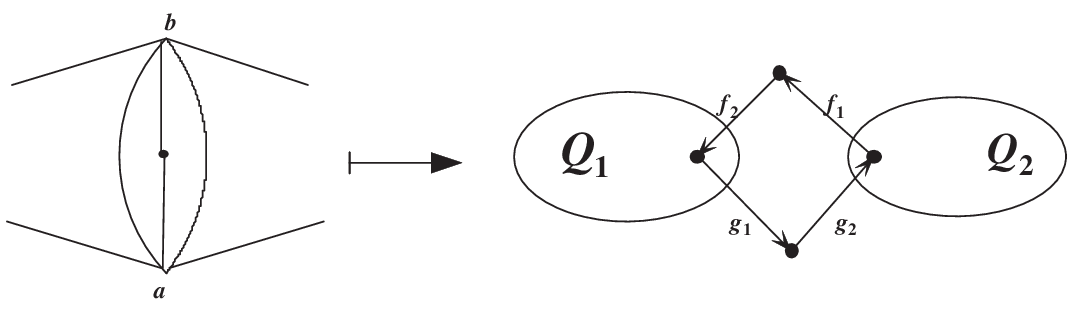} \\
{Figure 6}
\end{center}
\end{figure}

In this case,  we may write the two degenerate tagged edges as
$M_{a,a}^1$ and $M_{b,b}^1$.  According to Lemma 3.3 and Proposition
3.4,  we know that the two vertices $a$ and $b$ are not neighbors
and that $M_{a,b}$, $M_{b,a}\in \Delta$.

As in Figure 6, $Q_1$ and $Q_2$ are in $\mathcal{M}^A$, the
corresponding cluster-tilted algebra $H_\Delta$ is the quotient
algebra $kQ_\Delta/I$ with the relation ideal $I$ generated by the
set $\Sigma_{Q_1}\bigcup \Sigma_{Q_2}\bigcup
\{f_1f_2g_1,f_2g_1g_2,g_1g_2f_1,g_2f_1f_2\}$.

\vskip 0.2in

\textbf{Type 4}: {\it  $\Delta$ has exactly $t$ degenerate tagged
edges with $t\geq 3$.}

\vskip 0.1in

In this case, we may write them as $M_{a_1,a_1}^1$, $M_{a_2,a_2}^1$,
..., $M_{a_t,a_t}^1$ in the counterclockwise direction. For
convenience, we set $t+i\equiv i ( \ {\rm mod}\ t)$ and
$M_{a_i,a_{i+1}}=0$ if $a_i$ and $a_{i+1}$ are neighbors. According
to Lemma 3.3,  we know that the non-zero elements of $M_{a_1,a_2}$,
$M_{a_2,a_3}$, ..., $M_{a_{t-1},a_t}$, $M_{a_t,a_1}$ belong to
$\Delta$.  Moreover, $M_{a_1,a_1}^1$, $M_{a_2,a_2}^1$, ...,
$M_{a_t,a_t}^1$ divide the triangulation $\Delta$ to $t+1$ parts
$\Delta_1$, $\Delta_2$, ..., $\Delta_t$ and $\{M_{a_1,a_2},
M_{a_2,a_3},
...,M_{a_{t-1},a_t},M_{a_t,a_1},M_{a_1,a_1}^1,...,M_{a_t,a_t}^1\}$
as in Figure 7. Note that the corresponding subquivers $Q_1$, $Q_2$,
..., $Q_t$ of $Q_\Delta$ are all in $\mathcal{M}^A$.  we denote
$V_{M_{a_i,a_i}}$ by $V_i$ and  $V_{M_{a_i,a_{i+1}}}$ by $V_{i\,
i+1}$.

\begin{figure}[htbp]
\begin{center}
\includegraphics[width=12cm,height=6cm]{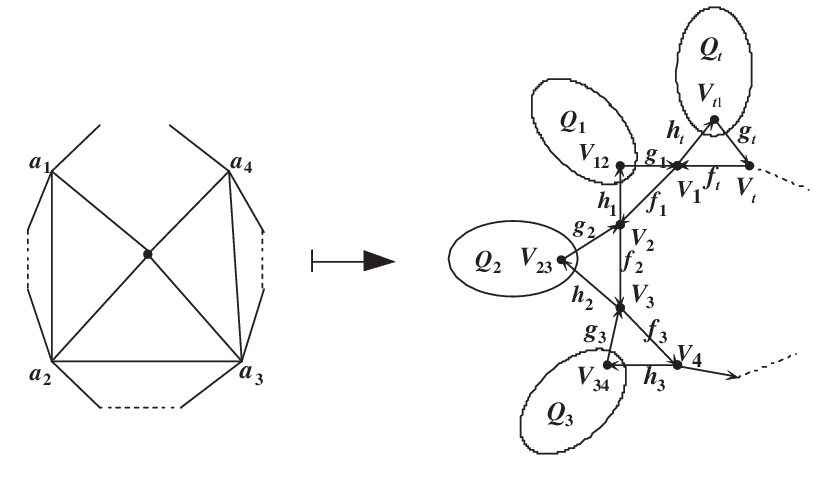} \\
{Figure 7}
\end{center}
\end{figure}

Note that $Q_i$ has no vertex if $a_i$ and $a_{i+1}$ are neighbors,
and that  $Q_i$ has only one vertex $V_{i\, i+1}$ if
$|\delta_{a_i,a_{i+1}}|=3$. The corresponding cluster-tilted algebra
$H_\Delta$ is the quotient algebra $kQ_\Delta/I$ with the relation
ideal $I$ generated by the set $(\bigcup_{i=1}^t
\Sigma_{Q_i})\bigcup \Sigma'\bigcup\Sigma''$, where $\Sigma'$ is the
set

$\{f_1g_1,f_2g_2, ..., f_tg_t, g_1h_1,g_2h_2, ..., g_th_t,
h_1f_1,h_2f_2, ..., h_tf_t\}$,

and $\Sigma''$ is the the union of the following sets

$\{ f_i f_{i+1} f_{i+2}...f_{i+t}, \ {\rm for} \ 1 \leq i \leq t, \
\ a_{i+t}\ {\rm is \ not \ a\ neighbor \ of} \ a_i,\}$

and $ \{ f_i f_{i+1} f_{i+2}...f_{i+(t-1)}, \ {\rm for} \ 1 \leq i
\leq t ,\ \  \ a_{i+t}\ {\rm is \ a\ neighbor \ of} \ a_i\}$.

\vskip 0.2in

We summarize the above discussion as the following main theorem
which has been proved by D. Vatne in [Vat] in different way.

\vskip 0.2in

{\bf Theorem 3.5.}  \ {\it  Let $n$ be an integer with $n\geq 4$ and
$Q$ be a quiver of a cluster-tilted algebra of type $D_n$. Then $Q$
must belong to one of type 1, type 2 ,type 3 or type 4.}
 \hfill$\Box$

\vskip 0.2in

\section { Isomorphism classes of cluster-tilted algebras of type $D_n$}

\vskip 0.2in

In this section, we follow the ideas of H.Torkildsen and generalize
the results for the cluster tilted algebras of type $A_n$ in [T] to
the type of $D_n$.

\vskip 0.2in

 Let $\mathbf{P}_n$ be the punctured polygon with $n\geq
5$. We write the boundary vertices as $\{ 1, 2, \cdots, n\}$ with
counterclockwise direction.

\vskip 0.2in

Let $H$ be the path algebra $kQ_n$ and $Q_n$ be the following
quiver:
$$\xymatrix{ & & & & & n\\
         1 \ar[r]& 2 \ar[r]& ...\ar[r] & n-3 \ar[r] & {n-2}
\ar[ru] \ar[rd]\\  & & & & & n-1 .}
$$
Let $\mathscr{C}_H$ be the cluster category of $H$,  let $T$ and
$T'$ be two basic tilting objects in $\mathscr{C}_H$. In the sequel,
we will give a criterion for when ${\rm End}_{\mathscr{C}_H}(T)^{\rm
op}$ and ${\rm End}_{\mathscr{C}_H}(T')^{\rm op}$ are isomorphic.

\vskip 0.2in

We now recall the Auslander-Reiten quiver of  the cluster category
$\mathscr{C}_H$, see [Hap]. It is a stable quiver built from $n$
copies of $Q_n$, we denote it by $\Gamma_{\mathscr{C}_H}$.  The
vertices of $\Gamma_{\mathscr{C}_H}$ are
$V(\mathscr{C}_H):=\mathbb{Z}_n \times \{ 1, 2, \cdots, n \}$.  The
arrows are $\begin{array}{rr}
(i,j)\rightarrow (i,l) \\
(i,l)\rightarrow (i+1,l) \end{array}\}$ whenever there is an arrow
$j\rightarrow l$ in $Q_n$.

\vskip 0.2in

Finally, the translation $\tau$ is defined by
$$ \tau (i,j) = \left\{
\begin{array}{ll}
(i-1,n-1), & {\rm if}\ i=0, \ j=n \ {\rm and}\ n \ {\rm is \ odd},\\
(i-1, n),  & {\rm if} \ i=0,\ j=n-1\  {\rm and}\ n \ {\rm is \ odd}, \\
(i-1,j), & {\rm otherwise}.
\end{array}\right.
$$

\vskip 0.2in

Let $\mathcal{C}$ be the category of tagged edges. By Lemma 2.2,
there is an equivalence of triangulated categories $\varphi:\
\mathcal{C} \rightarrow \mathscr{C}_H$. We denote by $T_{i,j}$
 the indecomposable object in
$\mathscr{C}_H$ corresponding to the vertex $(i,j)$ in
$\Gamma_{\mathscr{C}_H}$.  It is easy to see that $\varphi
(M_{i,i+2}) = T_{i,1}$, $\varphi (M_{i,i}^1) = T_{i,n-1}$ and
$\varphi(M_{i,i}^{-1}) = T_{i,n}$.

\vskip 0.2in

In section 3, we have defined an automorphism $\sigma$ in
$\mathcal{C}$ which also induces an automorphism in $\mathscr{C}_H$,
and  we still denote it by $\sigma$. It follows that $\sigma
(T_{i,j}) = T_{i,j}$, for $1\leq j \leq n-2$, $\sigma(T_{i,n}) =
T_{i,n-1}$ and $ \sigma(T_{i,n-1}) = T_{i,n}$.

\vskip 0.2in

Let $\mathcal{T}_n^H$ be the set of all basic tilting objects of
$\mathscr{C}_H$. We define an equivalence relation $"\sim "$ on
$\mathcal{T}_n^H$ as following,  $T \sim T'$ if and only if $T =
\tau^i T'$ or $T =\sigma\tau^j T'$ for some $i$ and $j$. We denote
by $\widetilde{\mathcal{T}_n^H}$ the set of equivalence classes of
$\mathcal{T}_n^H$. Note that the  bijection $\varphi : \mathcal{T}_n
\rightarrow \mathcal{T}_n^H$ also induces a bijection
$\widetilde{\varphi}:\ \widetilde{\mathcal{T}_n} \rightarrow
\widetilde{\mathcal{T}_n^H}$.

\vskip 0.2in

Now we can state our main result in this section.

\vskip 0.2in

{\bf Theorem 4.1.}\ {\it Let $T$ and $T'$ be basic tilting objects
in the cluster category $\mathscr{C}_H$ of type $D_n$ with $n\geq
5$, then the cluster-tilted algebras
$\Gamma=\mathrm{End}_{\mathscr{C}_H}(T)^{\rm op}$ and
$\Gamma'=\mathrm{End}_{\mathscr{C}_H}(T')^{\rm op}$ are isomorphic
if and only if $T \cong \tau^i T'$ or $T\cong \sigma \tau^j T'$ for
some integers $i$ and $j$.}

\vskip 0.2in

In order to  prove the theorem, we need some more lemmas.

\vskip 0.2in

{\bf Lemma 4.2.}\ {\it  Let $\Delta$ be an element in
$\widetilde{\mathcal{T}_n}$ and $Q_\Delta =
\widetilde{\gamma}(\Delta)$ be the corresponding quiver. If there is
a tagged edge $M$ close to the border, then the corresponding vertex
$V_M$ in  $Q_\Delta$ either is a source, a sink, or lies on a
cycle.}

\vskip 0.1in

{\bf Proof.}\  We only need to  consider the following eight cases
as shown in Figure 8.

\begin{figure}[htbp]
\begin{center}
\includegraphics[width=9cm,height=9cm]{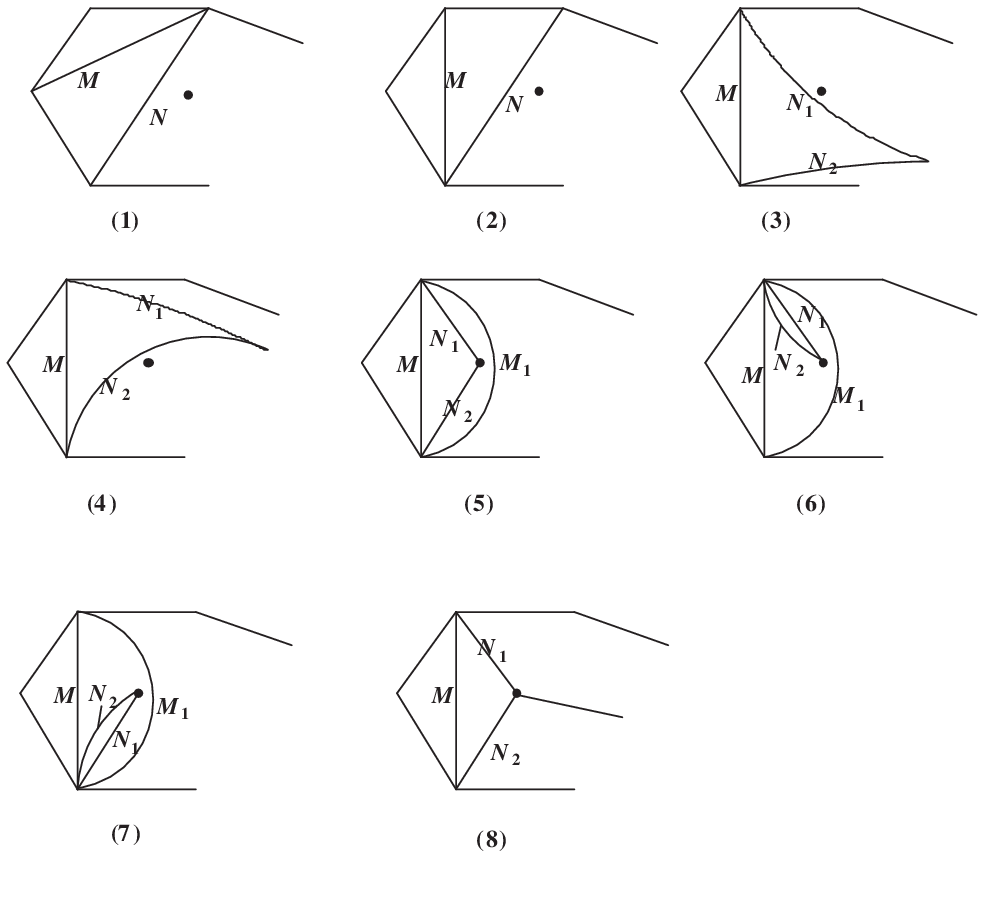} \\
{Figure 8}
\end{center}
\end{figure}

According to Lemma 3.2, it is easy to see that the vertex $V_M$ in
the first case is a source, and that the vertex $V_M$ in the second
case is a sink.  Finally,  vertex $V_M$ in the last six cases lies
on the following six cycles, respectively.

\vskip 0.2in
$\xymatrix@!=1.0pc{ V_M \ar[r] & V_{N_1} \ar[d] \\
                  & V_{N_2} \ar[lu]}$,
$\xymatrix@!=1.0pc{ V_M \ar[r] & V_{N_1} \ar[d] \\
                  & V_{N_2} \ar[lu]}$,
$\xymatrix@!=1.0pc{V_M \ar[r] & V_{N_1} \ar[d] \\
                  V_{N_2} \ar[u] & V_{M_1} \ar[l]}$,
$\xymatrix@!=1.0pc{V_M \ar[rd] & V_{N_1}\ar[l] \\
                   V_{N_2} \ar[u] & V_{M_1} \ar[l] \ar[u]}$,
$\xymatrix@!=1.0pc{V_M \ar[r] \ar[d] & V_{N_1}\ar[d] \\
                   V_{N_2} \ar[r] & V_{M_1} \ar[lu] }$ and
$\xymatrix@!=1.0pc{ V_M \ar[r] & V_{N_1} \ar[d] \\
                  & V_{N_2} \ar[lu]}$. $\Box$

\vskip 0.2in

\begin{figure}[htbp]
\begin{center}
\includegraphics[width=10cm,height=3cm]{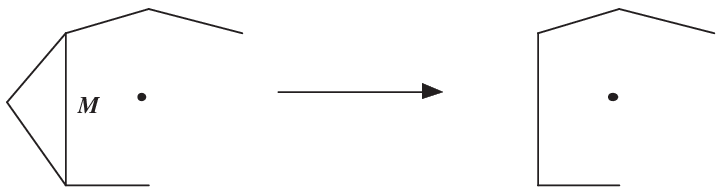} \\
{Figure 9}
\end{center}
\end{figure}

 Let $\Delta$ be a triangulation of the punctured polygon
$\mathbf{P}_n$ with $n\geq 5$ and let $Q_\Delta=
\widetilde{\gamma}(\Delta)$ be the corresponding quiver. Let $M$ be
a tagged edge in $\Delta$ and $V_M$ be the vertex in $Q_\Delta$
corresponding to $M$. We define a quotient triangulation $\Delta/M$
from $\Delta$ by factoring out $M$ and leaving all the other tagged
edges unchanged. Clearly, if $M$ is close to the border, $\Delta/M$
is a triangulation in the puncture polygon $\mathbf{P}_{n-1}$ as
Figure 9.

\vskip 0.2in

We denote by $Q_\Delta/V_M$ the quiver corresponding to $\Delta/M$.
According to Lemma 4.2, we have the following corollary.

\vskip 0.2in

{\bf Corollary 4.3.}\ {\it  Let $\Delta$ be a triangulation in
$\mathcal{T}_n$ with $n\geq 5$ and $Q_\Delta =
\widetilde{\gamma}(\Delta)$ be the corresponding quiver. Suppose
that $\Delta$ has a tagged edge $M$ which is  close to the border.
Let  $V_M$ be the vertex in $Q_{\Delta}$ corresponding to $M$.  Then
$Q_\Delta/V_M$ is a connected quiver in $\mathcal {M}_{n-1}^D$.}

\vskip 0.2in

{\bf Lemma 4.4.}\ {\it  Let $\Delta$ be a triangulation in
$\mathcal{T}_n$ with $n\geq 5$ and $Q_\Delta =
\widetilde{\gamma}(\Delta)$ be the corresponding quiver. Suppose
that $M$ is  a degenerate tagged edge in $\Delta$ and that $V_M$ is
the corresponding vertex to $M$ in $Q_{\Delta}$. Then $Q_\Delta/V_M$
is a connected quiver in $\mathcal {M}_{n-1}^A$.}

\vskip 0.1in

{\bf Proof.}\  According to Theorem 3.5, we only need to consider
the following four cases as shown in Figure 10.

\begin{figure}[htbp]
\begin{center}
\includegraphics[width=10cm,height=5cm]{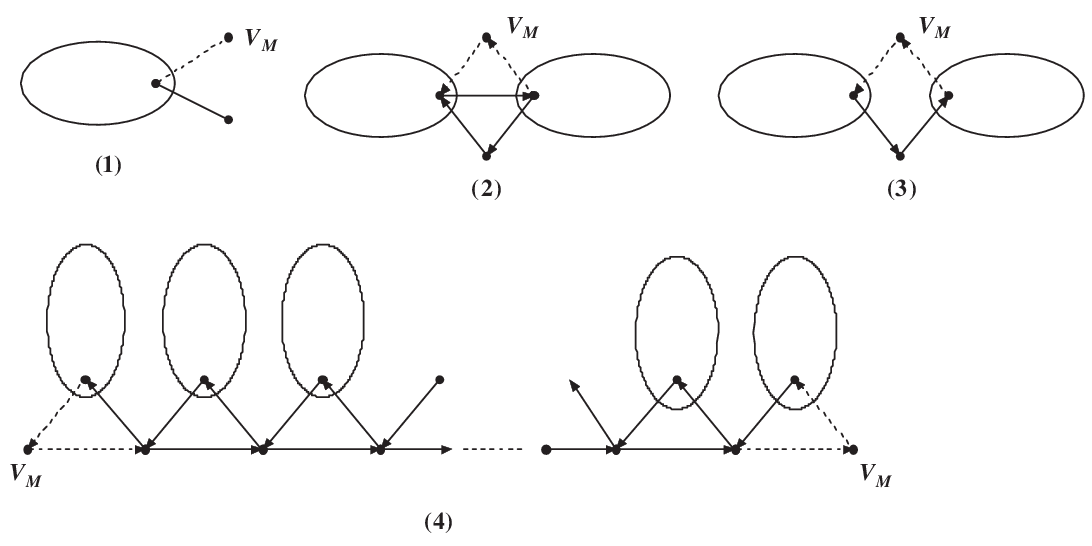} \\
{Figure 10}
\end{center}
\end{figure}

By [BV], our statement is true.  \hfill$\Box$

\vskip 0.2in

{\bf  Remark.} \ If $\Delta$ has a tagged edge $M$  which is
connected, that is, the corresponding vertex $V_M$ in $Q_{\Delta}$
is connected, by using Lemma 3.1 and Lemma 3.2,  we know that the
quiver $Q_\Delta/V_M$ is not connected.

\vskip 0.2in

{\bf Proposition 4.5.}\  {\it Let $\Delta$ be an element in
$\mathcal{T}_n$ with $n\geq 5$ and $Q_\Delta =
\widetilde{\gamma}(\Delta)$ be the corresponding quiver. Let $M$ be
a tagged edge in $\Delta$ and $V_M$ be the vertex in $Q_{\Delta}$
corresponding to $M$. Then}

\vskip 0.1in

(1). {\it  $Q_\Delta/V_M$ is connected and in $\mathcal{M}_{n-1}^D$
if and only if $M$ is close to the border.}

\vskip 0.1in

(2). {\it $Q_\Delta/V_M$ is connected and in $\mathcal{M}_{n-1}^A$
if and only if $M$ is degenerate.}

\vskip 0.1in

{\bf Proof.}\ It follows from Corollary 4.3, Lemma 4.4 and the above
Remark. \hfill$\Box$

\vskip 0.2in

{\bf Lemma 4.6.} \ {\it The function $\tilde{\gamma}:\
\widetilde{\mathcal{T}_5} \rightarrow \mathcal{M}_5^D$ is
injective.}

\begin{figure}[htbp]
\begin{center}
\includegraphics[width=15cm,height=8cm]{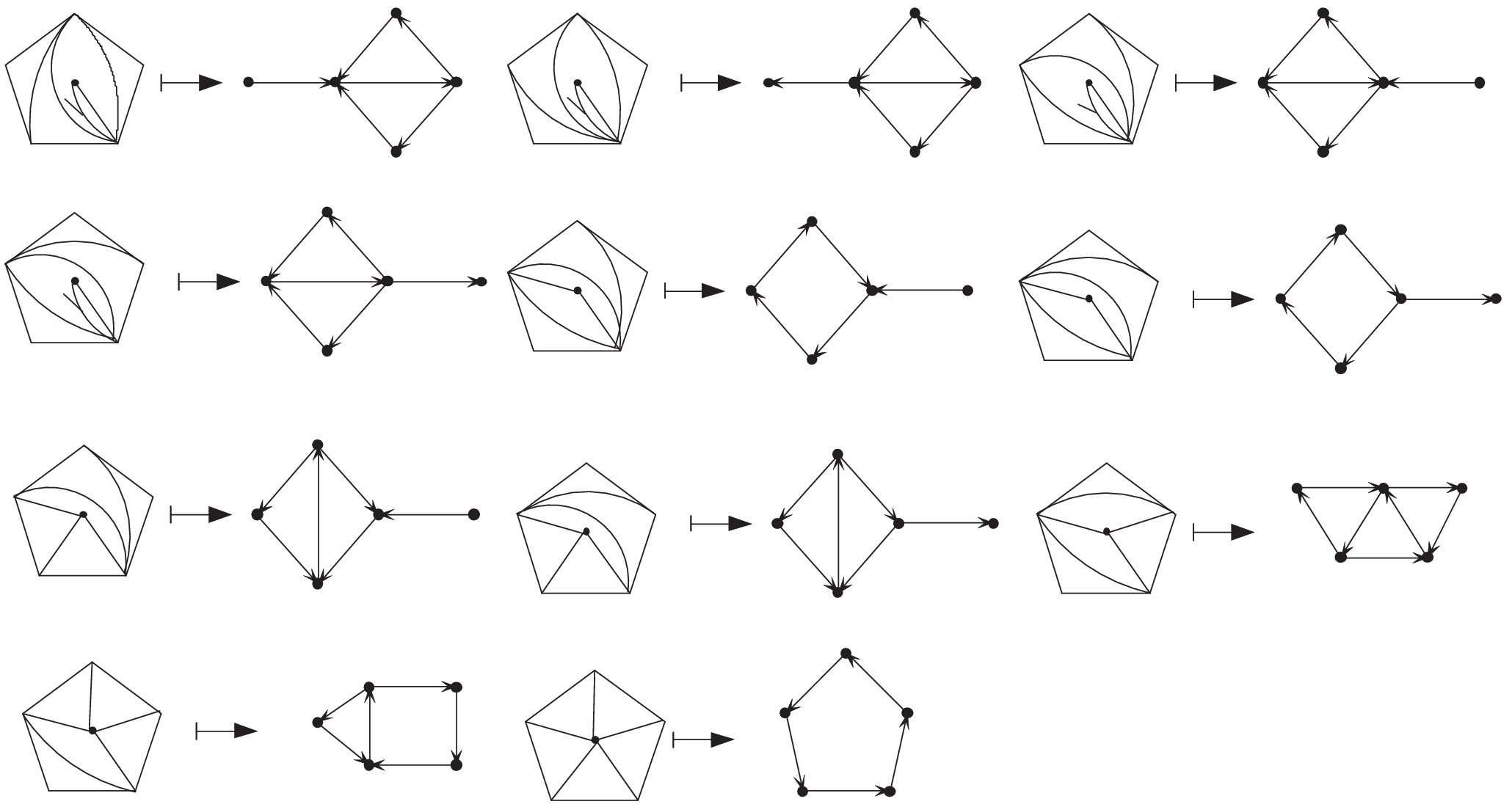} \\
{Figure 11}
\end{center}
\end{figure}

{\bf Proof.}\  There are only four types of elements in
$\widetilde{\mathcal {T}_5}$ by the discussion in section 3. First,
it is easy to see that the number of all the triangulations of Type
1 in $\widetilde{\mathcal {T}_5}$  are fifteen, and all of them are
mapped by $\tilde{\gamma}$ to non-isomorphic quivers. Next, we list
all the triangulations of Type 2, Type 3 and  Type 4 in
$\widetilde{\mathcal {T}_5}$, and their corresponding quivers as in
Figure 11. It is easy to see that all of them are mapped by
$\tilde{\gamma}$ to non-isomorphic quivers. In particular, we have
that the function $\tilde{\gamma}$: $\widetilde{\mathcal
{T}_5}\rightarrow \mathcal{M}_5$ is injective. \hfill$\Box$

\vskip 0.2in

Now, we can prove the result promised at the end of section 2.

\vskip 0.2in

{\bf Proposition 4.7.}\ {\it The function $\tilde{\gamma}:\
\widetilde{\mathcal{T}_n} \rightarrow \mathcal{M}_n^D$ is bijective
for all integer $n\geq 5$.}

\vskip 0.1in

{\bf Proof.}\  Note that $\tilde{\gamma}$ is always surjective, it
is sufficient to prove that it is injective. Suppose that
$Q_{\Delta}=\tilde{\gamma}(\Delta) = \tilde{\gamma}(\Delta')=
Q_{\Delta'}$ in $\mathcal{M}_n^D$. We will show that $\Delta
=\Delta'$ in $\widetilde{\mathcal{T}_n}$ by induction on $n$.

For $n=5$, by Lemma 4.6, we know that $\tilde{\gamma}:\ \widetilde{
\mathcal{T}_5} \rightarrow \mathcal{M}_5^D$ is injective.

Now suppose  that $\tilde{\gamma}:\ \widetilde{ \mathcal{T}_{n-1}}
\rightarrow \mathcal{M}_{n-1}^D$ is injective.

If $\Delta$ has no tagged edge close to the border, then every
tagged edge in $\Delta$ is degenerate, it is easy to see that
$Q_{\Delta}=Q_{\Delta'}$ is an oriented cycle with $n$ vertices.
 It forces that every tagged edge in $\Delta'$ is also degenerate,
 hence $\Delta =\Delta'$ in $\widetilde{\mathcal{T}_n}$.

 \begin{figure}[htbp]
\begin{center}
\includegraphics[width=13cm,height=10cm]{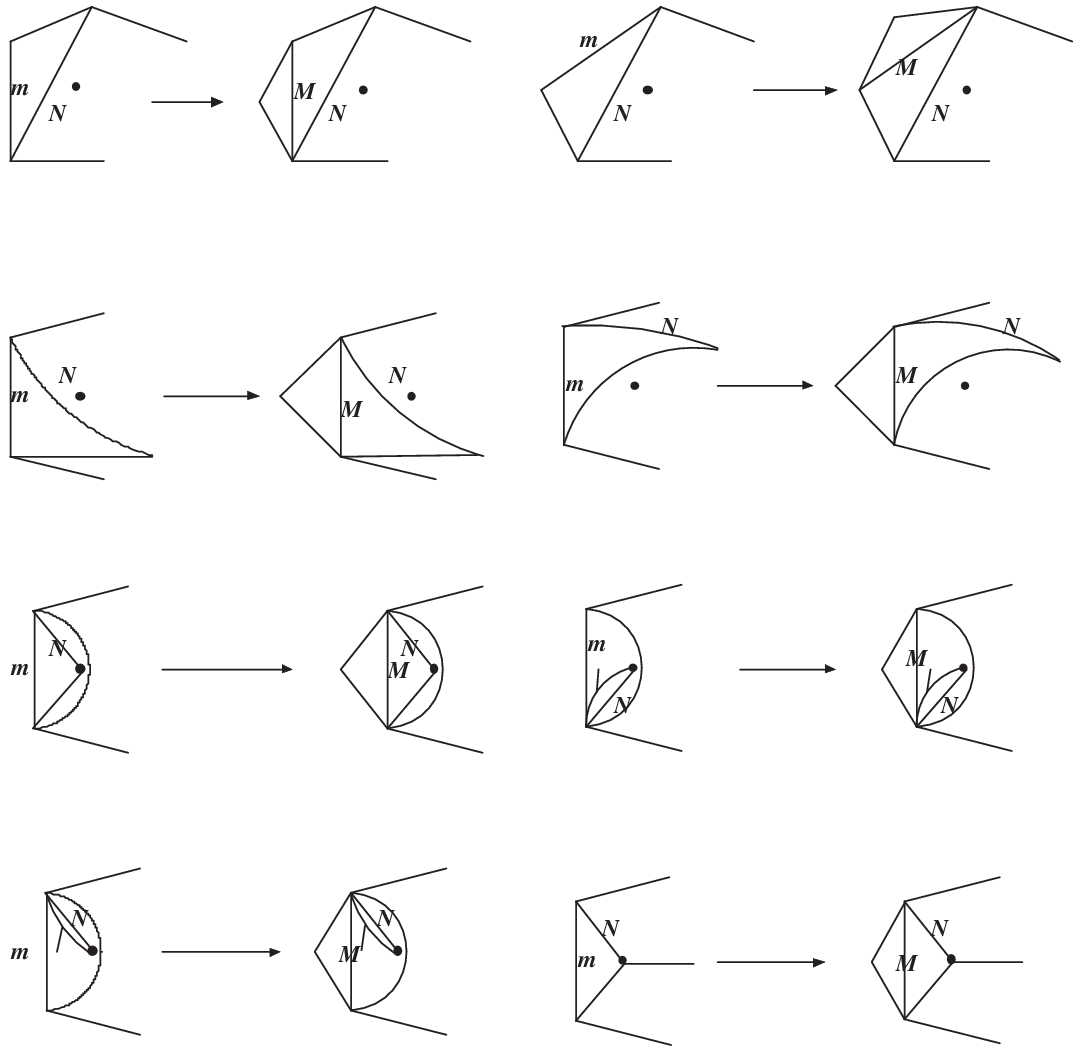} \\
{Figure 12}
\end{center}
\end{figure}

If $\Delta$ has a tagged edge $M$ close to the border, let $V_M$ be
the vertex in $Q_\Delta$ corresponding to $M$.  Note that
$Q_{\Delta}=Q_{\Delta'}$,  according to proposition 4.5,  we know
that $\Delta'$ must have a tagged edge $M'$ which is also close to
the border and that $M'$ is corresponding to $V_{M'}=V_M$ in
$Q_{\Delta}=Q_{\Delta'}$. According to Corollary 4.3, We know that
$Q_{\Delta/M} = Q/V_M = Q'/V_{M'} = Q_{\Delta'/M'}$ belongs to
$\mathcal{M}_{n-1}^D$,  and by the inductive hypothesis, We know
that $\Delta/M =\Delta'/M'$ in $\widetilde{\mathcal{T}_{n-1}}$.

\vskip 0.1in

We can obtain $\Delta$ and $\Delta'$ from $\Delta/M =\Delta'/M'$ by
extending the punctured polygon $\mathbf{P}_{n-1}$ at some border
edge. Fix a tagged edge $N$ in $\Delta$ such that $V_M$ and $V_N$
are adjacent in $Q_\Delta$. Let $N'$ be the tagged edge in $\Delta'$
corresponding to $V_N$.

As in Figure 12, there are at most eight ways to extend $\Delta/M$
such that the new tagged edge is adjacent to $N$. It is easy to see
that these extensions will be mapped by $\tilde{\gamma}$ to
non-isomorphic quivers. Also there are at most eight ways to extend
$\Delta'/M'$ such that the new tagged edge is adjacent to $N'$, and
all these extensions are mapped to non-isomorphic quivers and thus
$\Delta =\Delta'$ in $\widetilde{\mathcal{T}}_n$. \hfill$\Box$

\begin{figure}[htbp]
\begin{center}
\includegraphics[width=8cm,height=10cm]{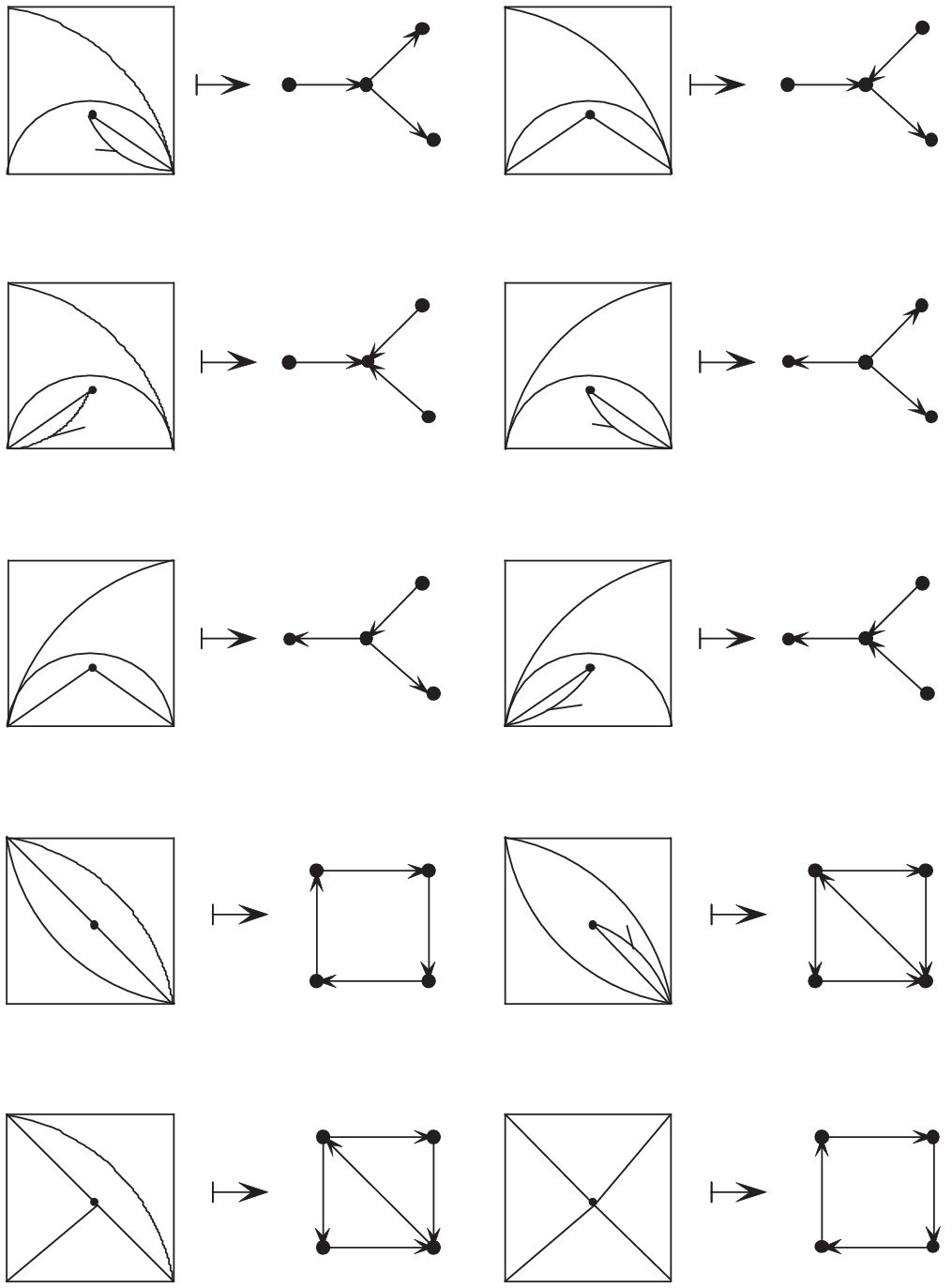} \\
{Figure 13}
\end{center}
\end{figure}

\vskip 0.2in

{\bf Remark.}\ Proposition 4.7 dose not hold for $D_4$ as showing in
Figure 13.

\vskip 0.1in

{\bf Proof of Theorem 4.1.}    Let $\Delta$ and $\Delta'$ be
 triangulations in $\mathcal {T}_n$, let  $T$ and $T'$ be the
 corresponding basic tilting objects  in $\mathscr{C}_H$.
 If $T \sim T'$, it is clear that
$\Gamma=\mathrm{End}_{\mathscr{C}_H}(T)^{\rm op}$ and
$\Gamma'=\mathrm{End}_{\mathscr{C}_H}(T')^{\rm op}$ are isomorphic.
Conversely, if $T\nsim T'$, then   $\Delta \nsim \Delta'$. By
Proposition 4.7, we know that the quiver $Q_\Delta$ is not
isomorphic to the quiver $Q_{\Delta'}$. It forces that
$\Gamma=\mathrm{End}_{\mathscr{C}_H}(T)^{\rm op}$ and
$\Gamma'=\mathrm{End}_{\mathscr{C}_H}(T')^{\rm op}$ are not
isomorphic, since every cluster-tilted algebra of type $D_n$ is
uniquely determined by its quiver.  This completes the proof.
\hfill$\Box$

\vskip 0.2in

{\bf Acknowledgements.}  After completing this work, the third
author was informed by Aslak Buan and Hermund Torkildsen that they
also proved the parallel results with this paper in [BT], he is
grateful to them for this.

\begin{description}

\item{[ARS]}\ M.Auslander, I.Reiten, S.O.Smal$\phi$, \ Representation
Theory of Artin Algebras. Cambridge Univ. Press, 1995.

\item{[ASS]}\ I.Assem, D.Simson, A.Skowronski, Elements of the
representation theory of associative algebras, Vol. 1, Cambridge
Univ. Press, 2006.

\item{[BMR1]}\ A. Buan, R. Marsh, I. Reiten,
               Cluster-tilted algebra.   Trans. Amer. Math. Soc.,
               359(1)(2007), 323-332.

\item{[BMR2]}\  A. Buan, R. Marsh, I. Reiten,
                 Cluster-tilted algebra of finite representation
                 type.   J. Algebra, 306(2)(2006), 412-431.

\item{[BMRRT]}\  A. Buan, R. Marsh, M. Reineke, I. Reiten, G.
                  Todorov,
           Tilting theory and cluster combinatorices,
               Adv.in Math., 204(2)(2006), 572-618.

\item{[BT]}\  A. Buan, H. Torkildsen, The number of elements in the mutation
         class of a quiver of type $D_n$,  arxiv:0812.2240, 2008.

\item{[BV]} \ A. Buan, D. Vatne,
                  Derived equivalence classification of
                 cluster-tilted algebra of type $A_n$.   J.
                 Algebra, 319(7)(2008), 2723-2738.

\item{[CCS1]}\  P. Caldero, F. Chapoton, R. Schiffler,
                   Quivers with relations arising from
                  clusters($A_n$ case).   Trans. Amer. Math. Soc.,
                   358(3)(2006), 1347-1364.

\item{[CCS2]}\  P. Caldero, F. Chapoton, R. Schiffler,
                    Quivers with relations and cluster-tilted
                      algebras.    Algebra. Represention. Theory,
                       9(2006), 359-376.

\item{[FZ1]} \ S. Fomin, A. Zelevinsky,
           Cluster algebra I: Foundation.   J. Amer. Math.
           Soc.,  15(2002), 497-529.

\item{[FZ2]} \   S. Fomin, A. Zelevinsky,
             Cluster algebra II: Finite type classification.
               Invention Math.,  154(2003), 63-121.

\item{[Hap]}\ D.Happel, Triangulated categories in the representation
theory of finite dimensional algebras.  Lecture Notes series 119.
Cambridge Univ. Press,1988.

\item{[Rin]}\ C.M.Ringel, Tame algebras and integral quadratic forms.
 Lecture Notes in Math. 1099. Springer Verlag, 1984.

\item{ [Sev]}\  A. Seven,
                   Recognizing cluster algebra of finite type.
                    Electron. J. combin., 14(1)(2007).

\item{[Sch]} \  R. Schiffler,
                  A geometric model for cluster category of type
                 $D_n$,   J. Algebra Comb., 27(1)(2008), 1-21.

\item{[Tor]}\  H. Torkildsen,
                    Counting cluster-tilted algebras of type
                      $A_n$,  International electronic Journal of
                      algebra,   4(2008), 149-158.

\item{[Vat]} \   D. Vatne,
                    The mutation class of $D_n$ quivers.
                   arxiv:0810.4789, 2008.

\end{description}

\end{document}